\documentstyle[12pt]{article}
\begin{document}
\begin{center}{\bf\LARGE  Quantization of function algebras on
semisimple orbits in ${\cal G}^*$}
\medskip
\newline Joseph Donin  and Steve Shnider, Bar Ilan University
\newline Dmitry Gurevich, Universit\'e de Valenciennes
\medskip
\newline{\em Preliminary version}
\end{center}
\baselineskip=20pt
\parskip=5pt
\parindent=18pt
\setcounter{page}{1}
\def\hsp{\hspace*{18pt}}
\newcommand{\QED}{\hspace{0.2in}\vrule width 6pt height 6pt depth 0pt 
\vspace{0.1in}} 
\newcommand{\Cal}{\cal}
\newcommand{\pf}{{\em Proof.} \hspace{0.2in}} 
\newcommand{\tag}[1]{\eqno(#1)}
\newtheorem{theorem}{Theorem}[section] 
\newtheorem{lemma}[theorem]{Lemma} 
\newtheorem{proposition}[theorem]{Proposition} 
\newtheorem{corollary}[theorem]{Corollary} 
\newtheorem{problem}[theorem]{Problem} 
\newtheorem{definition}[theorem]{Definition}
\def\t{\otimes}
\def \G{{\cal G}}
\def \A{{\cal A}}
\def \qi{q^{\frac{H_i}2}}
\def \qi-{q^{-\frac{H_i}2}}
\def \P{{\cal P}}
\def\E{\End(M_{\P})}
\def \N{{\cal N}}
\def \H{{\cal H}}
\def \F{{\cal F}}
\def \I{{\cal I}}
\def \O{{\cal O}}
\def \U{{\cal U}}
\def \L{{\cal L}}
\def\th{{\rm th}}
\def\C{{\bf C}}
\def\ch{K[[h]]}
\def\la{\lambda}
\def\End{{\rm End}}
\def\Ind{{\rm Ind}}
\def\Aut{{\rm Aut}}
\def\Ad{{\rm Ad}}
\def\dim{{\rm dim}}
\def\Ker{{\rm Ker}}
\def\Im{{\rm Im}}
\def\Ind{{\rm Ind}}

\begin{abstract}

In this paper we describe a multiparameter deformation of   
the function algebra of a  semisimple coadjoint orbit.
In the first section we use  the representation of the Lie algebra on a
generalized Verma module to quantize the Kirillov bracket on the family of
 semisimple coadjoint orbits of a given orbit type.
In the second section we extend this construction
to define a deformation in the category of representations of the quantized
enveloping algebra. In an earlier paper
 we used cohomological methods to prove the existence of a two
parameter family quantizing a compatible pair of Poisson brackets on 
any symmetric coadjoint orbit. This paper gives a more explicit algebraic
construction which includes more general orbit types and which we prove to be
flat in all parameters.
\end{abstract}

\begin{section}{Quantizing the Kirillov bracket}
Assume we have the following data: 
 a simple Lie algebra, $\G$, over the complex field, $\C$, with 
corresponding Lie group, $G$, and a Cartan subalgebra, $\H$, together with
a system of simple positive roots. Let $\G=\N^-\oplus \H\oplus
\N^+$, be the corresponding Cartan decomposition,
where $\N^+,\N^-$ are the nilpotent subalgebras which are,
respectively, the sum of positive root spaces and the sum of  the
negative root spaces. Denote by $\G^*$ and  $\H^*$  the $\C$ linear
duals of $\G$ and $\H$. Given $\la\in\G^*$, let
$\O_\la$ be the orbit of $\la$ under the coadjoint action of $G$
on $\G^*$. Finally, let $S(\G)$ be the symmetric algebra of
$\G$ considered as polynomial functions on $\G^*$ and $\I_\la$ the
ideal of polynomials vanishing on the orbit $\O_\la$. Then  $\F_\la=
S(\G)/\I_\la$ is the algebra of functions on $\O_\la$, with Poisson
structure given by the standard Kostant-Kirillov-Souriau (KKS) bracket 
\newline $\{f,g\}(\la)=\langle\la, [df_{\la},dg_{\la}]\rangle$.

Consider the $G$ orbits of  semisimple elements,
that is,  linear functions
conjugate under the coadjoint action to ones
 which vanishes on $N^+\oplus \N^-$.
In any such orbit we can pick as  origin an element, $\la$,
which is the trivial extension to $\G$ of a
functional on $\H$. The stabilizer of such
an element will be a subalgebra generated by  $\H$ and
a subset of the simple roots. Such a subalgebra is called a 
Levi subalgebra.  The set of all $\la$ with
stabilizer equal to a fixed Levi algebra $\L$
is parametrized by a subspace of $\H^*$ minus its
intersection with a family of coordinate hyperplanes. Denote
this parameter space by $\Lambda_\L$, and represent the
elements by $m$-tuples, $\la=(\la_1,\ldots, \la_m)$,
where the $\la_i$ are the values of $\la$ on a fixed subset 
of the set of simple positive roots and  all $\lambda_i$ are different from
zero. (In the case of $\G=sl(n,\C)$, this corresponds to
a Levi algebra which is block diagonal with $m+1$ blocks.)
The algebras we are interested in quantizing are the $\F_\la$
as $\la$ varies in $\Lambda_\L.$

Define the parabolic subalgebra, $\P=\L\oplus\N^+$.
Let  $M_\la=\Ind_{\P}^{\G}{\bf 1}_\la$ be the representation
of $\G$ induced from  the one dimensional representation of $\P$
defined by $\la.$  Let $\N^-_{\P}$ be the subalgebra 
of $\N^-$ which is a vector space complement to
$\P$ in $\G$.
As a vector spaces all the $M_\la$
can be identified with $M_{\P}:=U(\N^-_{\P})$ by the isomorphism
$$M_\la=U(\G)\otimes_{U(\P)}{\bf 1}_\la\cong U(\N^-_{\P})=M_{\P}.$$ 
Henceforth we use this  identification and denote  by $v_0$ the
element of $M_\la$ corresponding to $1\in U(\N^-_{\P})$.
The element $v_0$ will be the highest weight vector of $M_\la$ for
all $\la.$ Let $\phi_\la:\G\rightarrow \E$  be the Lie algebra homomorphism
defining the  representation $M_\la$, then
for all $X\in \G$ and all $\la\in \Lambda_{\L}$
$$\phi_\la(X)v_0=\la(X)v_0+\ldots\mbox{lower weight terms}.$$ 
We study the quantization relative to indeterminates,
$h$ ( the formal deformation parameter) and 
 $\la_i$ (parametrizing $\Lambda_\L$).
Let $\End(M_{\P})[\la,h]$ be the ring 
of polynomials  $h,\la_1,\ldots,\la_m$ with coefficients in $\End(M_{\P})$
and $T({\G})[\la, h]$ the ring of polynomials in $h,\la_1,\ldots,\la_m$
with coefficients in the tensor algebra of $\G$. The map 
$$h\phi_{\la/h}:\G\longrightarrow \E[\la,h]$$
extends to a homomorphism of $\C[\la, h]$ algebras:
$$ \phi_{\la,h}:T({\G})[\la,h]\longrightarrow \E[\la, h].$$
Consider $T({\G})[\la,h]$ as a graded algebra with elements of
$\G$ and the indeterminates $\la$ and $h$ considered as having degree
one. Define a grading  on the algebra $\End(M_{\P})[\la,h]$ with
the elements of $\End(M_{\P})$ having degree zero and $\la$ and $h$
having degree one. Then let $\A_{\la,h}$ be the subalgebra of
$\End(M_{\P})[\la,h]$  generated by the image of
$\phi_{\la,h}$. Since $\phi_{\la,h}$ is a map of graded algebras, the image
$\A_{\la,h}$ is a graded subalgebra.  As a torsion free  submodule of  the free 
$\C[h]$ module, $\E[\la,h]$, it is also a  free  $\C[h]$
module.  We will prove the following
\begin{theorem}
The graded algebra $\A_{\la,h}$ is a free  $\C[\la,h]$ module.
At  all points $\la=\la_0$ of $\Lambda_\L$ , $\A_{\la_0,h}$ is a quantization
of the Poisson algebra $\F_{\la_0}$ with KKS bracket.
\end{theorem}
\pf In the course of proving the second part of 
the theorem we will show that the quotient $A_{\la,h}/hA_{\la,h}$ is
isomorphic to a certain free $\C[\la]$ module, $N\t_\C \C[\la]$,
where $N$ is an induced representation of $\G$ defined
independently of $\la.$ Then
$$A_{\la,h}\cong A_{\la,h}/hA_{\la,h}\t_\C \C[h]\cong N\t_\C\C[\la]
\t_\C\C[h]\cong N\t_\C\C[\la,h],$$ 
which proves the first assertion of the theorem.

In order to prove the second statement we must prove the following two
statements: 
\begin{enumerate}
\item
For all $\la_0\in\Lambda_\L$ the algebra $\A_{\la_0,h}/h\A_{\la_0,h}$
is isomorphic to the algebra  $\F_{\la_0}$.
\item
Denoting the isomorphism in 1. by $\sigma$
and defining the Poisson bracket
of two cosets,
$[f],[g]$, in $\A_{\la_0,h}/h\A_{\la_0,h}$, by  
$$\{[f],[g]\}=[\frac1h(fg-gf)]\quad\mbox{we obtain}\quad\sigma\{[f],[g]\}
 =\{\sigma([f]),\sigma([g])\}.$$
\end{enumerate}

When $\la$ is evaluated at a fixed vector in
$\Lambda_\L$, the algebra $\A_{\la,h}$ is, by definition, generated 
as an $\C[h]$ subalgebra of $\E[h]$ by the elements
$\phi_{\la,h}(X)=h\phi_{\la/h}(X)$ for  $X\in\G$. 
To simplify notation, let $\phi:=\phi_{\la/h}$. The Lie algebra
homomorphism property 
$\phi(X)\phi(Y)-\phi(Y)\phi(X)=\phi([X,Y])$
implies that 
$$
(h\phi)(X)(h\phi)(Y)=(h\phi)(Y)(h\phi)(X)\quad
\mbox{mod}\quad h\cdot\Im(h\phi).$$
Therefore the algebra $\A_\la=\A_{\la,h}/h\A_{\la,h}$ is commutative and
$\phi_{\la,h}$ induces a map 
$$\psi:S(\G)\rightarrow \A_{\la}.$$ 
Let $X_i^{\pm}$ be a basis for $\N^{\pm}$ and $X^0_j$ a basis for
$\H$, and, as usual, for any multi-index $A=(a^-_i,a^0_j,a^+_i)$,
let $X^A=\Pi(X^-_i)^{a^-_i}\Pi(X^0_j)^{a^0_j}\Pi(X^+_i)^{a^+_i}$
and $|A|=\sum a^+_i +\sum a^0_j +\sum a^-_i.$
We use $X^{\tilde A}$ to denote  the symmetrization of the
monomials $X^A$, considered
as elements of $T(\G)$ and as a basis for $S(\G)$. By definition
$$\phi_{\la,h}( X^{\tilde A})=h^{|A|}(\phi(X))^{\tilde A},$$
where the latter expression is the symmetrization in $\E$ of the products of the linear factors 
$h\phi(X_i)$. 
Now suppose that $f(X)=\sum_A c_A(\la)
X^{\tilde A}\in \Ker\psi$,  then
there exists an element of $g(X)\in T(\G)[\la,h]$ 
 such that $\phi_{\la,h}(f(X))=h\phi_{\la,h}(g(X)).$
Applying the difference to the highest weight vector, we get
$$\phi_{\la,h}(f(X)-hg(X))v_0=0.$$ On the 
other hand, the coefficient of the $v_0$ 
in this expression depends only on the 
summands  which have weight zero relative to $\H$. 
Furthermore any term of weight zero
containing a pair of a negative root vector, $\phi(X^-_i),$
and a positive root vector, $ \phi(X^+_i)$, 
will contribute to the coefficient of $v_0$ only via commutators 
$\phi([X^-_i,X^+_i])$. For such terms the homogeneity in $h$ 
 exceeds the number of factors containing $\phi$. The same is true of
all the terms in $\phi_{\la,h}(hg(X))$.
This means that the   constant term in $h$  in the coefficient of $v_0$
in $\phi_{\la,h}( f(X)-hg(X))v_0$ is  
$$0=\sum_A c_A(\la(X))^{ A}=f(\la(X)).$$
But identifying $S(\G)$ with the polynomial functions
on $\G^*$, $f(\la(X))$ is the  $f(X)(\la)$,
 the  evaluation of the function
$f(X)$ at the point $\la\in\G^*$. This proves that any
$f(X)\in \Ker\psi$ vanishes at the point $\la$. 
 To prove that
$f(X)$ vanishes on the entire orbit, 
we use the invariance of $\Ker\psi$ 
under the adjoint action of $G$. For any 
$g\in G$ we have $f(X)(g\la)=
(g^{-1} f(X))(\la)=0, $ which completes the proof that
 $\Ker \psi\subset \I_{\la}$.

The inclusion just proved shows that $\psi$ induces an epimorphism
\begin{equation}\label{2}
\pi_\la: \A_\la\cong S(\G)/\Ker\psi\rightarrow S(\G)/\I_{\la}=\F_{\la}.
\end{equation}

To prove that this is an isomorphism, we use the fact that
$\F_\la$ decomposes into a direct sum of 
finite dimensional irreducible representations of $\G$. 
On the other hand, being  a free $\C[h]$ module,
 $A_{\la,h}\cong N\t_\C\C[h]$ and as
$\C$ vector spaces, $A_\lambda=A_{\la,h}/hA_{\la,h}\cong N$.
 Evaluation at the generic value of $h$ gives an 
isomorphism of $N$ with the image of
the representation $\phi_{\la/h}$, where, for any $\la\in \Lambda_\L$,
$\la/h$ is a generic $\C$ valued  weight vector in $\Lambda_\L$. 
Let
 $$N=\oplus n_\mu(\la/h) V_\mu \quad\mbox{and}\quad
\F_\la=\oplus  \ell_\mu(\la) V_\mu.$$
The $\G$ equivariance of the  epimorphism $\pi_\la$ implies that
$n_\mu(\la/h)\geq \ell_\mu(\la).$

To establish the reverse inequality and thus to show that
$\pi_\la$ is an isomorphism, we construct a $\G$ equivariant
imbedding of $N$ in $\F_\la$.  The construction depends on the
 following lemma. 
For details we refer the
reader to  \cite{D}
\begin{lemma}
 Let $M_{\la}^*$  be the restricted dual of $M_{\la}$.
For generic $\la$, $M_{\la}^*$ is a generalized Verma module
for the parabolic algebra $\L\oplus \N_\P^-$
generated as a principal $U(\N^+_\P)$ module by the lowest weight vector
$v_0^*.$ 
\end{lemma}

Let $\Ind_\L^\G{\bf 1}$ be the representation of $\G$ induced
from the trivial representation of $\L$ with basis vector ${\bf 1}$.
 Then
${\bf 1}\mapsto v_0
\otimes v_0^*$ induces a $\G$ equivariant morphism
$$\tau_{\la}:\Ind_\L^\G {\bf 1}\stackrel{\cong}{\longrightarrow}
M_{\la}\otimes_{\C} M_{\la}^*.$$  
We have the following isomorphisms 
and identities,
$$\Ind_\L^\G{\bf 1}\cong U(\N^+_\P)U(\N^-_\P){\bf 1}, \quad
M_\la\cong U(\N^-_\P)v_0,
\quad U(\N^-_\P)v_0^*=0,$$
and for generic $\la$
$$M^*_\la\cong U(\N^+_\P)v_0^*.$$ 
A simple inductive argument using the filtration by weights in
$M_\la$ and $M_\la^*$ shows
that for generic  $\la$, $\tau_{\la}$ is an isomorphism.

Dualizing this isomorphism and restricting to 
$\G$ finite elements in $\End(M_{\la})$ relative to the conjugation
action,   we get an isomorphism
$$\tau_\la^*: (M_{\la}\otimes
M_{\la}^*)^{*\G-finite}\stackrel{\cong}{\longrightarrow} (\Ind_\L^\G{\bf
1})^{*\G-finite}.$$

An element of $M_{\la}\otimes M_{\la}^*$ defines an element of
finite rank in $\End(M_{\la})$ in the obvious way. We can define a
nondegenerate pairing between  $A\in\End(M_{\la})$
and $B\in M_{\la}\otimes M_{\la}^*$ 
by $\langle A,B\rangle=trace (A\circ B)$. This induces an imbedding $j$ of
$\End(M_{\la})$ into $(M_{\la}\otimes M_{\la}^*)^*$
whose restriction to the $\G$ finite elements will also
be denoted by $j$. 

For the generic value $\la/h$ the composition  $\tau_{\la/h}\circ j$
 defines an imbedding
$$\tau_\la\circ j:(\End(M_{\la/h}))^{\G- finite}\hookrightarrow
(\Ind_\L^\G{\bf 1})^{*\G-finite}.$$

The representation $\Ind_\L^\G{\bf 1}$ is naturally 
identified with point distributions on the orbit $\O_\la$ supported at the point
$\la$. Dualizing, we get an identification of
$(\Ind_\L^\G{\bf 1})^{*\G-finite}$ and  $\F_{\la}$. 

The restriction of the imbedding $\tau_{\la/h}\circ j$ to the
subset $N\subset
\End(M_{\la/h})$  composed with the
identification of $(\Ind_\L^\G{\bf 1})^{*\G-finite}$ and  $\F_{\la}$
defines  a $\G$ equivariant imbedding of $N$ in $\F_{\la}$. This proves
the reverse inequality relating multiplicities, 
$n_\mu(\la/h)\leq \ell_\mu(\la).$ 
Therefore we have equality 
$$n_\mu(\la/h)= \ell_\mu(\la),$$ 
which implies that $\pi_\la$ is an
isomorphism and $\Ker\psi_\la=\I_\la,$ concluding the proof 
that $A_{\la,h}/hA_{\la,h}\cong \F_\la.$ 

 When we consider 
$\la$ as an indeterminate, the isomorphism
$$\F_\la\cong(\Ind_\L^\G{\bf 1})^{*\G-finite}\t_\C\C[\la]$$
also proves that $A_{\la,h}/hA_{\la,h}$ is a free $\C[\la]$ module.
 
Finally, in relation to the Poisson structures, it is enough to
check condition 2. on linear  elements in $A_{\la, h}$ and $\F_\la$.
For $X,Y\in\G$ considered as functions on $ \O_\la$ their Poisson
bracket  is the function defined by $[X,Y]$, the bracket in $\G$.
On the other hand the bracket of the corresponding elements in
$A_{\la, h}/hA_{\la, h}$ is 
$$\frac1h((h\phi)(X)(h\phi)(Y)-(h\phi)(Y)(h\phi)(X))=h\phi([X,Y]),$$
as required.\QED
\end{section}

\begin{section}{Quantization in the category of $U_q(\G)$ modules}

Given $\G$ as above, let $U_q(\G)$ be the quantized 
enveloping algebra  with generators, $E_i,F_i,K^\pm_i,\, 1\leq i\leq n,$ 
satisfying the standard relations given by Lusztig, see \cite{Lu}. 
We adopt the widely used convention  of using $q^{H_i}$ to represent
$K_i.$ 
Let $\U^-$ be the subalgebra generated by the $F_i$, $\U^0$, the
(commutative) subalgebra generated by the $q^{H_i}$ which, as an
abelian group relative to multiplication, is 
isomorphic to the lattice $\Pi$ generated by the simple roots,  and $\U^+$ the
subalgebra generated by the $E_i$. 
The algebra $\U$ decomposes into a product $\U=\U^-\U^0\U^+$. For any 
Levi subalgebra $\L$, the quantized universal enveloping algebra
$U_q(\L)$ is  naturally imbedded in $U_q(\G)$.  Let $\U_\L$ denote
the  image of this imbedding. It is generated by $E_i, F_i$ for
$i$ in a subset of the simple roots, and all $K_i^{\pm}$. 
Define the subalgebra of
$\U$ associated to the parabolic subalgebra $\P\supset \L$ by
$$\U_\P:=U_\L \U^+.$$

 For any weight vector $\lambda=\sum\lambda_i \omega_i\in \Lambda_\L$,
let ${\bf 1}_{q,\lambda}$ be the one
dimensional representation of $\U_\P$ with basis vector $v_0$,
which is the trivial representation when restricted to $\U_\L^-$ and to $\U^+$
and is defined  for the remaining generators of $\U_\P$, that is, 
those $q^H\in \U^0$,  by 
$$q^H v_0=q^{\langle\lambda, H\rangle} v_0.$$ 
Inducing this representation up to $\U$ 
gives the quantized generalized Verma module: 
$$M_{q,\lambda}=\U\t_{\U_\P}{\bf 1}_{q,\lambda}\quad 
\phi_{q,\lambda}(u) (w\t v_0)=uw\t v_0, \quad u,w\in \U.$$
In this construction we may consider
$q$  and the components of the vector $\lambda$ 
 as complex parameters,  or we may
let $q=e^t$ and consider all constructions over the formal 
power series in $t$ and $\lambda_it$. 
As a vector space over $\C$ or a module over $\C[[t,\lambda_it]]$,
$M_{q,\lambda}$ is isomorphic to  $M_\L:=\U^-/\U_\L^-$.

The next step is to find 
a subrepresentation of $\End(M_{q, \lambda})$ which is
 a deformation of the adjoint representation of $\G$.
Previously this followed trivially from the fact that the subspace  
$\G\subset U(\G)$ is  $ad(U(\G))$ invariant.  
In the present context we rely on the work of Joseph and Letzter, 
\cite{JL}, on the structure of the $ad(\U)$ finite part of 
$\U$, denoted $F(\U)$. 

First of all, the subspace $ad(\U)q^{\alpha}$ is  finite dimensional
if and only if $\alpha=-4\mu$ for a nonnegative integral weight $\mu.$
Let $R$ be the set of such elements of the lattice, $\Pi,$ and 
 $F(\mu)= ad(\U)q^{-4\mu}$for $\mu\in R$. Then 
$$F(\U)=\oplus_{\mu \in R} F(\mu)
\quad\mbox{and }\quad F(\mu)\cong \End(V_{q,\mu}),$$
where $V_{q,\mu}$ is the representation of $\U$ deforming 
the finite dimensional representation of $\G$ with highest weight
$\mu.$ 

Let $\G_q$ denote the deformed
adjoint representation,  which, by elementary deformation theory,
is unique up to equivalence. 
The representation $\G_q$ occurs infinitely often in $F(\U)$ 
since each $\End(V_{q,\mu})$ contains a copy of $\G_q$.
However, see \cite{JL}, for any finite dimensional representation,  
$V_{q,\mu}$, the  isotypical component, $F(\U)_{V_{q,\mu}}$, 
is a finitely generated module over the center, with multiplicity equal to
 the dimension of the zero weight space in the nonquantized 
representation, $V_{\mu}.$ Thus, up to multiplication by scalar
operators, there are finitely many copies of
$\G_q$ in  the image of $F(\U)$ in  $\End(M_{q,\lambda})$. 

We fix a $\G_q\subset F(\U)$ which specializes at $q=1$ to
the standard imbedding of $\G$ in $U(\G)$, where specialization at $q=1$,
is carried out most simply by replacing the generators
$q^{\pm H_i}$ with $(q^{\pm H_i}-1)/(q^{d_i}-q^{-d_i})$ so
that the Hopf algebra structure is defined over $\C[q, q^{-1}].$
Although one could construct the quantization  we are seeking in terms of the 
parameter $q$ it is convenient to revert to the setting of formal deformations,
setting $q=e^t$ and considering $t$ and the components $\lambda_i$ as
indeterminates. We continue to use the symbol $q$ but for the
remainder of this section it will be understood as a formal
power series in $t$, and  $\G_q$ will be a free $\C[[t]]$ module.
The tensor algebra $T(\G_q)$ is defined relative to the tensor product over
$\C[[t]]$. 

Let $\End(M)[[t,\lambda/h,h]]=\End(M_{q,\lambda/h})\t\C[[h]]$
considered as a $\C[[t,\lambda_i, h]]$ module.
From the representation of $F(\U)$ on the quantized Verma module we get a map
$$\phi=h\phi_{q,\lambda/h}:\G_q[[\lambda,h]]\rightarrow
\End(M)[[t,\lambda/h, h]].$$
This extends to a $\C[[t,\lambda/h,h ]]$ module map of the
tensor algebra $$\phi:T(\G_q)[[\lambda,h]]\rightarrow
\End(M)[[t,\lambda/h, h]].$$
We can define  a grading on $T(\G_q)[[ \lambda,h]]$
by setting the elements of $\G_q$, and the variables $\lambda_i$ and $h$
to be  of  degree one and $t$ to be  of degree zero. Similarly, we have
a grading on $\End(M)[[t,\lambda/h, h]]$ defined in the  same way on the
parameters and  with elements of $\End(M)$ considered as degree zero.
Then $\phi$ is a morphism of graded modules.

Define the subalgebra  $A_{t,\lambda,h}:=\mbox{Im}\,\phi.$
At $t=0$ we recover the construction of section 1
quantizing the Kirillov bracket. Now, the
dependence on $t$ introduces a new deformation parameter 
which quantizes some quadratic bracket  arising 
from the $U_q(\G)$ module structure.

\begin{theorem}
The algebra $A_{t,\lambda,h}$ is a $U_q(\G)$ module, and is free as a
$\C[[t,\lambda, h]]$ module.
Multiplication is $U_q(\G)$ equivariant, 
\end{theorem}
\pf
The definition of the action of   $u\in\U$ on
$A\in\End(M)[[t,\lambda,h]]$,
$$ad(u)\cdot A=\sum \phi(u_{(1)})\circ A\circ \phi(S(u_{(2)})
\quad\mbox{implies}\quad \sum ad(u_{(1)})(A)ad(u_{(2)})(B)=ad(u)(AB).$$
Since $\phi(\G_q)$ is a $\U$ submodule so is the span of all products
of elements in $\phi(\G_q)$, but this is precisely $A_{t,\lambda,h}.$
Moreover the latter formula is precisely the definition of equivariance
of multiplication. To prove that $A_{t,\lambda,h }$ is free we use the fact 
that the quotient module, 
$A_{t,\lambda,h}/tA_{t,\lambda,h}\cong A_{\lambda,h}$,
that is, the algebra defined in section 1, is a free
$\C[\lambda,h]$ module. Hence the same is true of the extension
to formal power series. Moreover $A_{t,\lambda, h}$ is free
 as a $\C[[t]]$ module
since it is a submodule of the free $\C[[t]]$ module, 
$\End(M)[[t,\lambda/h,h]]$.  Hence
$$A\cong A/tA\hat\t\C[[t]\cong N\hat \t \C[[\lambda, h]]\hat\t \C[[t]]\cong
 N\hat \t \C[[t,\lambda, h]].$$
{\bf Remarks 1.} The algebra $A_{t,\lambda, h}$ is graded as a quotient of 
$T(\G_q)[[\lambda,h]]$ by a graded ideal.\newline
{\bf 2.} The bracket defined by the commutator modulo the ideal generated by
$t$ and $h$ is a  linear combination of the Kirillov bracket and a second
bracket which is some kind of $R$-matrix bracket on the orbit. We have 
verified that for $\G=sl(2)$, the second bracket is the $R$ matrix bracket
generated by $E\wedge F$, see \cite{DS},  using
the results of Lyubashenko and Sudberry on quantum Lie algebras, 
see \cite{LS}. We 
believe that this is the case for the general, $\G$ and  symmetric coadjoint
orbit. It is interesting question to 
determine the form of the second bracket for
the general semisimple orbit, since it is known that the naive definition of
an $R$-matrix bracket is not in general compatible with the Kirillov bracket.
\end{section}

\end{document}